\documentclass[12pt,leqno]{amsart}
\usepackage{amsfonts,amssymb,mathrsfs,graphicx,amsmath, color}
\usepackage[utf8]{inputenc}   % LaTeX, comprends les accents !
\usepackage[T1]{fontenc}      % Police contenant les caract??res français
\usepackage{geometry}         % Définir les marges
\usepackage[francais]{babel}  % Placez ici une liste de langues, la
\newtheorem{theo}{Th\'{e}or\`{e}me}
\newtheorem{defi}{D\'{e}finition}
\newtheorem{propo}{Proposition}

\newtheorem{lemm}{Lemme}

\newcommand{\K}{\mathbb{K}}
\newcommand{\C}{\mathbb{C}}
\newcommand{\g}{\frak{g}}

\newcommand{\R}{\mathbb{R}}

\newcommand{\ds} {\displaystyle}

\newcommand{\h}{\frak{h}}

\newcommand\pp{\mathcal{P}}

\newcommand\dd{\noindent{\it D\'emonstration. }}

\newcommand\be{\begin{enumerate}}
\newcommand\ee{\end{enumerate}}

\title{Alg\' ebro\"{\i}des de Lie. Alg\`ebres de Lie-Rinehart}
\author{Elisabeth Remm}
\date{}
\address{ Universit\'e de Haute Alsace, Laboratoire IRIMAS-Math\'ematiques, 6 rue des FR\`eres Lumi\`ere. F.68093 Mulhouse}
\email{elisabeth.remm@uha.fr}

\begin{document}

\maketitle

\noindent{INTRODUCTION.} 
Les algèbres de Lie-Rinehart constituent un modèle algébrique des différentes structures que possèdent le fibré tangent à une variété différentielle, à savoir l'ensemble de ses sections est une algèbre de Lie, c'est aussi un module sur l'algèbre associative commutative des fonctions différentiables sur la variété et il agit sur cette même algèbre comme des dérivations (c'est la définition même d'un champ de vecteurs). Géométriquement, cette structure sur le fibré tangent est appelé un algébroïde de Lie. La généralisation de cette notion à un fibré vectoriel est dûe à J. Pradines \cite{Pra}.  Certaines structures géométriques sur une variété permettent de définir des algébroïdes de Lie associés. Dans ce travail, on construit de tels algébroïdes pour des structures de contact ou plus généralement pour certains systèmes de Pfaff de contact. On s'intéresse ensuite à l'aspect algébrique en étudiant des algèbres de Lie-Rinehart, tout d'abord en petite dimension, puis lorsque l'algèbre de Lie associée est l'algèbre des dérivations de l'algèbre associative commutative donnée par la structure de Lie-Rinehart. On montre enfin que certaines algèbres de Courant permettent de définir des algèbres de Lie-Rinehart qui vues comme une algèbre ordinaire (définie par une seule multiplication) est une algèbre de Leibniz gauche. 

\tableofcontents
\section{Alg\'ebro\"{\i}de de Lie: D\'efinition et rappels}

\subsection{Le fibr\' e tangent \`a une vari\'et\'e diff\'erentielle}
Soit $V$ une vari\' et\' e diff\' erentielle de dimension $n$ et soit $x_0$ un point de $V$. Soit $\mathcal{C}^\infty (V)$ l'alg\`ebre associative commutative des fonctions diff\' erentiables sur $V$. 
\begin{defi}
Un vecteur tangent en $x_0$ est une application
$$D: \mathcal{C}^\infty (V) \rightarrow \R$$
v\' erifiant
\begin{enumerate}
\item $D(f+g)=D(f)+D(g), \ \ \forall f,g \in \mathcal{C}^\infty (V),$
\item$ D(fg)=D(f)g(x_0)+f(x_0)D(g), \ \ \forall f,g \in \mathcal{C}^\infty (V),$
\item $D(f)=0$ si $f$ est constante.
\end{enumerate}
\end{defi}
L'ensemble de ces applications appel\'ees aussi bien  d\' erivations en $x_0$ que vecteurs tangents en $x_0$ forme une espace vectoriel r\' eel not\' e $T_{x_0}(V)$ appel\' e l'espace vectoriel tangent \`a $V$ en $x_0$. Dans un ouvert de coordonn\' ees $(x^1,\cdots,x^n)$ de $V$, les applications $(\frac{\delta}{\delta x^i})_{x_0}$ d\' efinies par
$$\displaystyle \left(\frac{\delta}{\delta x^i}\right)_{x_0}(f)=\left(\frac{\delta f}{\delta x^i}\right)_{(x^1_0\cdots,x^n_0)}$$
est une base de $T_{x_0}(V)$.
Rappelons \' egalement qu'un fibr\' e vectoriel de base $V$ et de fibre type $\R^p$ est un triplet $(E,\pi,V)$ o\`u $E$ est une vari\' et\' e de dimension $n+p$, $\pi$ une application diff\' erentiable surjective de $E$ sur $V$ tel que pour tout $x \in V$, la fibre $E_x=\pi^{-1}(x)$ est un espace vectoriel de dimension $p$. Le fibr\' e tangent \`a $V$ est le fibr\' e vectoriel ainsi construit: $E= \cup_{x \in V} T_x V$, l'application $\pi$ est la projection naturelle $\pi: w \in T_x(V) \rightarrow x$. On note le fibr\' e tangent
$$(T(V),\pi,V).$$
ou plus simplement $T(V)$.  Sa fibre type est donc $\R^n$.

Un champ de vecteurs $X$ sur $V$ est la donn\' ee pour chaque point $x \in V$ d'un vecteur $X_x \in T_x(V)$.  Il se comporte comme une application
$$X: \mathcal{C}^{\infty}(V) \rightarrow \mathcal{C}^{\infty}(V)$$
d\' efinie par $$X(f)(x)=X_x(f)$$
pour toute application diff\' erentiable $f$ sur $V$, 
On note g\' en\' eralement par $\mathcal{X} (V)$ l'ensemble des champs de vecteurs sur $V$.
Il est clair que $\mathcal{X} (V)$  est un espace vectoriel r\' eel (de dimension infini). C'est aussi un $\mathcal{C}^{\infty}(V)$-module de type fini, ceci signifiant que $\mathcal{X} (V)$ est engendr\' e (dans chaque ouvert de coordonn\' ees) par un nombre fini de champs de vecteurs $\displaystyle (\frac{\delta}{\delta x^i})$ les composantes \' etant des fonctions. On \' ecrira donc dans un ouvert de coordonn\' eees:
$$X= \displaystyle \sum_i X_i \frac{\delta}{\delta x^i}.$$
\begin{propo}
L'espace vectoriel $\mathcal{X} (V)$ est muni d'une structure d'alg\`ebre de Lie. Le crochet de Lie des champs de vecteurs $X$ et $Y$ est donn\' e par
$$[X,Y](f)=X(Y(f))-Y(X(f))$$
pour tout $f \in \mathcal{C}^{\infty}(V).$
Localement, ce crochet s'exprime ainsi
$$[X,Y]=\left[\sum_i X_i \frac{\delta}{\delta x^i},\sum_j Y_j \frac{\delta}{\delta x^j}\right]=\sum_{i,j}X_i  \frac{\delta Y_j}{\delta x^i}\frac{\delta}{\delta x^j} -Y_j  \frac{\delta X_i}{\delta x^j}\frac{\delta}{\delta x^i}.$$
\end{propo}

\subsection{L'alg\' ebro\"{\i}de de Lie "fibr\' e tangent"}
La description ci-dessus du fibr\' e tangent \`a $V$ nous montre
qu'il y a interaction entre l'alg\`ebre associative commutative $\mathcal{C}^{\infty}(V)$ et l'alg\`ebre de Lie $\mathcal{X}(V).$ Plus pr\' ecis\' ement
\begin{theo}
\begin{enumerate}
\item L'alg\`ebre de Lie $\mathcal{X}(V)$ est un $\mathcal{C}^{\infty}(V)$-module (projectif), c'est-\`a-dire, il existe une application
$$(f,X) \in \mathcal{C}^{\infty}(V) \times \mathcal{X}(V) \rightarrow fX \in \mathcal{X}(V)$$
v\' erifiant
$$[fX,gY]=fg[X,Y]+fX(g)Y-gY(f)X$$
pour tout $f,g \in \mathcal{C}^{\infty}(V)$ et $X,Y \in \mathcal{X}(V)$.
\item  $\mathcal{X}(V)$ est une alg\`ebre de Lie de d\' erivations de $\mathcal{C}^{\infty}(V)$: pour tout $X \in \mathcal{X}(V)$
$$X: \mathcal{C}^{\infty}(V) \rightarrow \mathcal{C}^{\infty}(V)$$
v\'erifiant
$$X(fg)=X(f)g+fX(g)$$
pour tout $f,g \in \mathcal{C}^{\infty}(V)$.
\item Les deux structures, celle d'alg\`ebre associative sur $\mathcal{C}^{\infty}(V)$ et celle d'alg\`ebre de Lie sur  $\mathcal{X}(V)$ sont li\'ees par une relation de Leibniz:
$$[X, fY ] = f[X, Y ] +X(f) Y.$$
\end{enumerate}
\end{theo}

\subsection{D\' efinition d'un alg\'ebro\"{\i}de de Lie}

\begin{defi}
On appelle alg\' ebro\"{\i}de de Lie un triplet $(E, [.,.], \rho)$ o\`u
\begin{enumerate}
\item $E$ est un fibr\' e vectoriel au dessus d'une vari\' et\' e $V$:\ $E\stackrel{\pi}\rightarrow V$
\item L'espace des sections $\Gamma(E)$ du fibr\' e $E$ est muni d'un crochet de Lie $[.,.]$ (une section est une application diff\' erentiable $s: V \rightarrow E$ telle que $\pi \circ s= Id$),
\item il existe un morphisme de fibr\' e vectoriel  $\rho:E \rightarrow T(V)$ qui induit un morphisme d'alg\`ebre de Lie  $\rho :\Gamma(E) \rightarrow \mathcal{X}(V)$ c'est-\`a-dire
$$\rho([X,Y]=[\rho(X),\rho(Y)]$$
pour tout $X,Y \in \Gamma(E)$
\item l'application $\rho$, appel\'ee ancre, et le crochet sur $\Gamma(E)$ v\' erifient la r\`egle de Leibniz:
$$[X,fY]=(\rho(X)(f))Y+f[X,Y]$$
pour tout $X,Y \in \Gamma(E)$ et $f \in \mathcal{C}^\infty (V).$
\end{enumerate}
\end{defi}
Notons que $\Gamma(E)$ est aussi un $\mathcal{C}^\infty (V)$-module via l'application $$(f \in \mathcal{C}^\infty (V), X \in \Gamma(E) )\rightarrow fX \in \Gamma(E)$$
avec $(fX)(x)=f(x)X(x)$ pour tout $x \in V$ et agit via l'application $\rho$ comme alg\`ebre de d\'erivations de $ \mathcal{C}^\infty (V)$.

\medskip

\noindent{\bf Exemples}
\begin{enumerate}
\item Le fibr\' e tangent est bien un alg\' ebro\"{\i}de de Lie, l'application ancre $\rho$ \' etant l'identit\' e.
\item Toute alg\`ebre de Lie est un alg\' ebro\"{\i}de de Lie. On consid\`ere dans ce cas que la vari\' et\' e $V$ est r\' eduite \`a un point
\item Vari\' et\' e de Poisson et alg\' ebro\"{\i}de de Lie.
Soit $V$ une vari\' et\' e diff\' erentielle et $\mathcal{C}^\infty (V)$ l'alg\`ebre des fonctions diff\' erentiables dur $V$. On sait que cette alg\`ebre est associative et commutative. On dit que $V$ est une vari\' et\' e de Poisson si l'on peut munir l'alg\`ebre $\mathcal{C}^\infty (V)$ d'un crochet de Lie $[.,.]$ (et donc $\mathcal{C}^\infty (V)$ est une alg\`ebre associative commutative et une alg\`ebre de Lie), ce crochet de Lie v\' erifiant la condition de Leibniz
$$[f,gh]=g[f,h]+[f,g]h$$
pour tout $f,g,h \in \mathcal{C}^\infty (V)$.
Ce crochet n'est pas \`a priori d\' efini \`a partir du crochet de Lie de $\mathcal{X}(V)$. Toutefois, si nous consid\' erons l'application
$$\Theta_f : \mathcal{C}^\infty (V) \rightarrow \mathcal{C}^\infty (V)$$
donn\' ee par
$$\Theta_f(g)=[f,g]$$
cette application, d'apr\`es la r\`egle de Leibniz, v\' erifie
$$\Theta_f (gh)=[f,gh]=g[f,h]+[f,h]g=g\Theta_f(h)+\theta_f(g)h.$$
C'est donc une d\' erivation de $\mathcal{C}^\infty (V)$ qui est donc associ\' ee \`a un champ de vecteurs de $V$. On peut n\' eanmoins associer \`a une vari\' et\' e de Poisson un alg\' ebro\"{\i}de de Lie en consid\' erant comme fibr\' e vectoriel sur $V$ le fibr\' e cotangent $T^*(V)$ dont chaque fibre $T_x^*(V)$ est l'espace dual \`a l'espace tangent $T_x(V)$. Une section de ce fibr\' e est une forme diff\' erentielle sur $V$. L'application ancre 
$$\rho : T^*(V) \rightarrow T(V)$$ est donn\' ee en coordonn\' ees locales par
$$\rho(\sum \alpha_i dx_i)= \sum \rho_{ij} \alpha_i \frac{\delta}{\delta x_j}$$
avec $\rho_{ij}=[x_i,x_j]$ o\`u $(x_1,\cdots,x_n)$ est un syst\`eme de coordonn\' ees locales de $V$ autour du point $x \in V$. On d\' efinit le crochet de Lie sur $T^* (V)$ de la fa\c{c}on suivante:
$$\{\alpha_1,\alpha_2\}=i(\rho(\alpha_1))d\alpha_2-i(\rho(\alpha_2))d\alpha_1+d(\Lambda(\alpha_1 \wedge \alpha_2)$$
o\`u $i(X)d\alpha (Y)=d\alpha (X,Y)$ et $\Lambda$ d\' esigne le tenseur de Poisson dont l'expression en coordonn\' ees locales est
$$\Lambda= \sum \rho_{ij}\frac{\delta}{\delta x_i} \wedge \frac{\delta}{\delta x_j}.$$
 
L'alg\' ebro\"{\i}de de Lie que nous venons de d\' efinir s'appelle l'alg\' ebro\"{\i}de de Lie de la vari\' et\' e de Poisson $V$. R\' eciproquement, \' etant donn\' ee une structure d'alg\' ebro\"{\i}de de Lie sur $T^*(V)$, si l'application $\Lambda$ associ\' ee \`a l'application ancre $\rho$ est antisym\' etrique, alors on peut d\' efinir sur $V$ une structure de vari\' et\' e de Poisson dont l'alg\' ebro\"{\i}de de Lie correspondant est celui qui est donn\' e.
\end{enumerate}

\section{Alg\' ebro\"{\i}des de Lie de contact}
\subsection{Alg\' ebro\"{\i}de de Lie associ\'e \`a une forme de contact}

 Soit $V$ une vari\' et\' e diff\' erentielle de dimension impaire $2p+1$. On dit que $V$ est une vari\' et\' e de contact s'il existe sur cette vari\' et\' e une forme de Pfaff $\omega$ (une forme diff\' erentielle de degr\' e $1$) non nulle en tout point de $V$ v\' erifiant l'identit\' e
$$\omega \wedge (d\omega)^p \neq 0.$$
Soit $Ker( \omega) = \cup_{x \in V} \ker \omega(x)$ la distribution d\' efinie par $\omega$. Cette distribution est non int\' egrable, c'est-\`a-dire n'est pas le fibr\' e tangent \`a une sous-vari\' et\' e de dimension $2p$ de $V$. Rappelons, que pour toute forme de contact sur $V$, la dimension des plus grandes sous-distributions r\' eguli\`eres (de dimension constante en tout point de $V$) de  $Ker (\omega)$ qui soient int\' egrables est inf\' erieure ou \' egale \`a $p$. Ainsi l'ensemble des sections du sous-fibr\' e $Ker (\omega)$ de $T(V)$ n'est pas une alg\`ebre de Lie pour le crochet des champs de vecteurs.
Comme $\omega$ est une forme de contact sur $V$, il existe un champ de vecteurs $R$ sur $V$, non nul en tout point, d\' etermin\' e par le syst\`eme
$$\omega(x)(R(x))=1, d\omega(R,X)=0$$
pour tout $x \in V$ et tout champ de vecteurs $X$ sur $V$.
Ce champ de vecteurs $R$ est appel\' e le champ de Reeb.  Ainsi, en tout point $x$ de $V$, nous avons la d\' ecomposition
$$T_x(V)=\R\{R(x)\}\oplus \ker \omega(x).$$
On d\' efinit ainsi un sous fibr\' e $R \rightarrow V$ de $T(V)$ dont la fibre au dessus du point $x$ est engendr\'ee par le vecteur $R(x)$.
Mais  \`a ce sous-fibr\' e  de $T(V)$ est associ\' e le fibr\' e quotient $T(V)/R$ d\' efini ponctuellement par $(T(V)/R)(x)=T_x(V)/R(x)$. Cet espace vectoriel est donc isomorphe \`a $\ker \omega(x)$. Notons par $E(\omega) \rightarrow V$ ce fibr\' e. Si $p$ d\' esigne la projection de $T(V)$ sur $E(\omega)$, alors le crochet de Lie de deux champs de vecteurs $\overline{X}=p(X), \ \overline{Y}=p(Y)$ de $E(\omega)$ s'\' ecrit
$$[\overline{X},\overline{Y}]=\overline{[X,Y]}$$
et l'ensemble des sections du fibr\' e $E(\omega)$ est une alg\`ebre de Lie. Consid\' erons \`a pr\' esent l'ensemble des fonctions diff\' erentiables sur $V$ v\' erifiant l'\' equation diff\' erentielle
$$R(f)=0.$$
Il est clair qu'elles constituent une alg\`ebre associative commutative et le couple $(L,A)$ o\`u $L$ est l'alg\`ebre de Lie des sections de $E(\omega)$ et $A$ l'alg\`ebre des fonctions v\' erifiant $R(f)=0$ est un alg\' ebro\"{\i}de de Lie associ\' ee \`a la vari\' et\' e de contact $(V,\omega)$.

Dans un ouvert de Darboux de $V$ (o\`u plus simplement si $V=\R^{2p+1}$, la forme $\omega$ s'\' ecrit
$$\omega=dx_1+x_2dx_3+\cdots+x_{2p}dx_{2p+1}$$
le champ de Reeb
$$R=\frac{\delta}{\delta x_1}.$$
Tout champ de vecteurs sur $E(\omega)$ peut s'\' ecrire
$$\overline{X}=\sum_{i=2}^{2p+1} X_i(x_2,\cdots,x_{2p+1})\frac{\delta}{\delta x_i}$$
et l'alg\`ebre associative $A$ est repr\' esent\' ee par les fonctions de $2p$ variables $f(x_2,\cdots,x_{2p+1}).$ On a bien
$$\overline{X}f=\sum_{i=2}^{2p+1} X_i(x_2,\cdots,x_{2p+1})\frac{\delta f}{\delta x_i}$$
et 
$$f\overline{X}=\sum_{i=2}^{2p+1} f X_i(x_2,\cdots,x_{2p+1})\frac{\delta }{\delta x_i}.$$
On v\' erifie que
$$
\begin{array}{lll}
\medskip
[\overline{X},f\overline{Y}]&=&\sum_{i,j} X_i\frac{\delta fY_j}{\delta x_i}\frac{\delta }{\delta x_j}-fY_j\frac{\delta X_i}{\delta x_j}\frac{\delta }{\delta x_i}\\
\medskip
&=&\sum_{i,j}f X_i\frac{\delta Y_j}{\delta x_i}\frac{\delta }{\delta x_j}+Y_jX_i\frac{\delta f}{\delta x_i}\frac{\delta }{\delta x_j}-fY_j\frac{\delta X_i}{\delta x_j}\frac{\delta }{\delta x_i}\\ 
\medskip
&=&f[\overline{X},\overline{Y}]+\overline{X}(f)\overline{Y}.
\end{array}
$$

Prenons comme exemple caract\' eristique de vari\' et\' e de contact le groupe de Lie de Heisenberg.  Il est repr\' esent\' e par le groupe des matrices
$$
\begin{pmatrix}
	1&x & z\\
	0&1&y\\
	0&0&1
\end{pmatrix}
$$
La forme de Pfaff $\omega=dz+xdy$ est une forme de contact et $R=\frac{\delta}{\delta z}$ est le champ de Reeb associ\' e. 

\subsection{Alg\' ebro\"{\i}de de Lie associ\' e \`a un syst\`eme de contact}

Un syst\`eme de Pfaff de rang $p$ sur une vari\' et\' e diff\' erentiable $V$ est un sous-fibr\' e vectoriel $E$ de rang $p$ du fibr\' e cotangent $T^* V$. Nous avons not\' e par $\chi(V)$ le module des champs de vecteurs sur $V$. Si $\Lambda^1(V)$ est l'espace des formes de Pfaff sur $V$ c'est-\`a-dire des $1$-formes ext\' erieures sur $\chi(V)$,  alors un syst\`eme de Pfaff peut \^etre vu comme un sous-module de $\Lambda^1(V)$ tel que pour tout point $x \in V$, 
$\{\alpha(x), \ \alpha \in E\}$ soit un sous-espace vectoriel de dimension $p$ de $T^*_x(V)$. On en d\' eduit que pour tout $x \in V$, il existe un ouvert $U$ de $V$ contenant $x$ tel que la restriction \`a $U$ du syst\`eme de Pfaff soit d\' efinie par $p$ formes de Pfaff, $\alpha_1,\cdots, \alpha_p$ lin\' eairement ind\' ependantes en tout point de $U$.

La notion duale de celle de syst\`eme de Pfaff est la notion de syst\`eme diff\' erentiel. Un syst\`eme diff\' erentiel $\mathcal{S}$ de rang $q$ sur $V$ est un sous-fibr\' e vectoriel du fibr\' e tangent $T(V)$. On peut aussi le d\' efinir comme un sous-module de $\chi(M)$ tel que pour tout point $x \in M$, $\{X(x), \ X \in \mathcal{S}\}$ est un sous-espace vectoriel de dimension $q$ de $T_x(M)$. Dans ce cas aussi, pour tout $x \in V$, il existe un ouvert  de $V$ contenant $x$ tel que la restriction de $\mathcal{S}$ \`a $U$ soit engendr\' ee par $q$ champs $X_1,\cdots,X_q$ lin\' eairement ind\' ependants en tout point de $U$. Cette dualit\' e entre syst\`eme de Pfaff et syst\`eme diff\' erentiel se r\' esume ainsi: si $\mathcal{S}$ est un syst\`eme diff\' erentiel sur $V$ de dimension $q$, alors son orthogonal $\mathcal{S}^\perp=\{\alpha \in \Lambda^1(V), \ \alpha(X)=0 \ \forall X \in \mathcal{S}\}$ est un syst\`eme de Pfaff de dimension $p=n-q$ o\`u $n=\dim V$. 

Un syst\`eme diff\' erentiel $\mathcal{S}$ de dimension $q$ sur $V$ est dit int\' egrable s'il est le fibr\' e tangent \`a une sous-vari\' et\' e $W$ de dimension $q$ de $V$. Il est donc muni naturellement d'une structure d'alg\' ebro\"{\i}de de Lie, l'application ancre \' etant donn\' ee par l'inclusion $i: W \hookrightarrow
V.$ 
Un syst\`eme de Pfaff $\mathcal{P}$ de rang $p$ sur $V$ est int\' egrable s'il existe une sous vari\' et\' e $W$ de $V$ de dimension $n-p$ et une immersion injective $i: W \hookrightarrow V$ telle que pour toute forme de Pfaff $\alpha \in \mathcal{P}$, on ait $i^*(\alpha)=0.$ Ceci est \' equivalent \`a dire que quelles que soient les formes de Pfaff $\alpha_1,\cdots,\alpha_p$ de $\mathcal{P}$, on ait pout tout $\alpha\in \mathcal{P}$
$$d\alpha \wedge \alpha_1 \wedge \cdots \wedge \alpha_p =0.$$

%%%%%%%%%%%%%%%%%%%%%
Consid\'erons \`a pr\' esent un syst\`eme de Pfaff de rang $p$ non int\' egrable. Soit $\{\alpha_1,\cdots,\alpha_n\}$ une base locale de formes de Pfaff sur $V$ telle que $\{\alpha_1,\cdots,\alpha_p\}$ soit une base du syst\`eme de Pfaff $\pp$. La notion de classe (au sens d'Elie Cartan) permet de mesurer la non int\' egrabilit\' e d'un tel syst\`eme.  Elle est d\' efinie en chaque point $x \in V$ comme suit: 
soit $F_x$ le sous-espace de $T^*_x(V)$ engendr\' e par $\{\alpha_{p+1}(x),\cdots,\alpha_n(x)\}$ et pour $i=1,\cdots,p$ soit $S_i(x)$ le plus petit sous-espace de $F_x$ tel que $d\alpha_i(x) \in \Lambda^2F_x$. Soit $S_x$ la somme vectorielle des sous-espaces $S_i(x)$ pour $i=1,\cdots,p$. Si $s_x=\dim S_x$, alors la classe de $\pp$ au point $x$ est l'entier
$$c_x(\mathcal{P})= s_x+ p$$
o\`u $p$ est le rang du syst\`eme de Pfaff $\pp$. Ainsi, si $\pp$ est int\' egrable, on a $c_x(\pp)=p$ en tout point. La non int\' egrabilit\' e se traduit donc par $c_x(\pp)-p \neq 0$ en certain point. 
\begin{defi}
Le syst\`eme de Pfaff $\pp$ sur $V$ est dit de contact si sa classe est constante et maximale en tout point, c'est-\`a-dire
$$c_x(\pp)=n$$
pour tout $x$ de $V$.
\end{defi}
Par exemple, si $p=1$ et $n=2k+1$, le syst\`eme de Pfaff est r\' eduit \`a une seule forme de Pfaff et ce syst\`eme est de contact si et seulement si cette forme est de contact. Mais, lorsque $p \geq 2$, il n'existe pas, contrairement \`a $p=1$, de mod\`ele local. Localement, une forme de forme de contact s'\' ecrit, d'apr\`es le th\' eor\`eme de Darboux:
$$\alpha=dx_1+x_2dx_3+\cdots+x_{2k}dx_{2k+1}$$
Mais ce type de r\' esultat n'existe pas d\`es que $q \geq 2$. Par exemple, dans cette \' etude fondamentale \cite{Cartan}, Elie Cartan d\' etermine tous les mod\`eles locaux des syst\`emes de Pfaff de classe constante dans $\R^n$ pour $n \leq 5$.
  
  Consid\' erons le sous-espace  vectoriel $S_x+\pp_x$     de $T^*_x(V)$ est donc de dimension $c_x(\pp)$. Il est appel\' e le {\bf syst\`eme caract\' eristique au point $x$ de $\pp$} et est not\' e $\mathcal{C}_x(\pp).$
 \begin{propo}
 Soient $X_{p+1},\cdots,X_n$ des champs de vecteurs locaux lin\' eairement ind\' ependants en tout point tels que
 $$\alpha_i(X_j)=0, \ \ i=1,\cdots,p,\ \ j=p+1,\cdots,n.$$
  Alors les formes de Pfaff
  $$\alpha_1,\cdots,\alpha_p, \imath(X_i)d\alpha_j$$
  engendrent  $\mathcal{C}_x(\pp)$ en tout point, o\`u $ \imath(X_i)d\alpha_j$ d\' esigne la forme de Pfaff d\' efinie par  $\imath(X_i)d\alpha_j(Y)= d\alpha_j(X_i,Y).$
  \end{propo}
 L'espace caract\' eristique au point $x \in V$ est l'orthogonal  $\mathcal{C}_x^\perp(\pp)$ de  $\mathcal{C}_x(\pp)$. Si on suppose que $\pp$ est de classe constante, on d\' efinit ainsi un syst\`eme diff\' erentiel int\' egrable. Mais ce syst\`eme est trivial si on suppose que $\pp$ est un syst\`eme de contact. 
 
 \noindent{\bf Remarque.} On peut \' egalement mesurer la non int\' egrabilit\' e du syst\`eme $\pp$ de la fa\c{c}on suivante: pour chacune des formes $\alpha_i$ de la base locale choisie de $\pp$, on a:
 $$\ds d\alpha_i= \sum_{j=1}^p \mu_i^j \wedge \alpha_j + \sum_{j,k=p+1}^n C_{j,k}^i \alpha_j \wedge \alpha _k$$
 et le syst\`eme est int\' egrable si et seulement si 
 $$\partial \alpha_i= \sum_{j,k=p+1}^n C_{j,k}^i \alpha_j \wedge \alpha _k=0$$
 pour $i=1,\cdots,p.$ Notons \' egalement que les formes $\mu_i^j $ d\' efinissent une connexion sur le fibr\' e vectoriel $\pp$ dont la courbure est donn\' ee par
 $$\Omega_i^j=d\mu_i^j+\sum_k\mu_i^k \wedge \mu_k^j.$$
 Nous \' etudierons plus tard ces connexions adapt\' ees au syst\`eme de Pfaff. 
 
 Consid\' erons \`a pr\' esent un syst\`eme de Pfaff $\pp$ de rang $p$ de classe maximum $n$ sur $V$ o\`u $n=\dim V$. D\`es que $n \geq 2$, un tel syst\`eme n'est pas int\' egrable, mais on peut s'int\' eresser aux vari\' et\' es int\' egrales de dimension plus petite que $n-p$, l'espace tangent en un point $y \in V$ sera donc un sous-espace vectoriel propre de $\pp_y$. Rappelons le r\' esultat
 \begin{theo}\cite{Go-CRAS} Soit $\mathcal{Q}$ un sous-fibr\' e de $\pp$ de rang $q$ d\' efinissant un syst\`eme de Pfaff int\' egrable sur $V$.Alors
 $$q \leq \ds \frac{p(n-p)}{p+1}.$$
 \end{theo}
 Comme nous savons, d'apr\`es \cite{Cartan} qu'il existe une infinit\' e de mod\`eles locaux de syst\`eme de Pfaff de contact de rang donn\' e, nous allons nous int\' eresser ici aux syst\`emes de Pfaff de contact ayant des vari\' et\' es int\' egrales de dimension maximale $q=  \ds \frac{p(n-p)}{p+1}.$
 
 \begin{defi}
 Soit $V$ une vari\' et\' e diff\' erentielle de dimension $n$ et $\pp$ un syst\`eme de Pfaff sur $V$ de rang $p$ et de classe maximale $n$. On dit que $\pp$ est un $p$-syst\`eme de contact s'il existe une sous-vari\' et\' e $W$ de $V$   int\' egrale \`a $\pp$ de dimension "maximale" $q= \ds \frac{p(n-p)}{p+1}$
c'est-\`a-dire si pour tout point $x \in V$, $T_xW$ est un sous-espace vectoriel de dimension  $q$ de $\pp_x$.
\end{defi}
Comme cons\' equence imm\' ediate, la dimension $n$ de la vari\' et\' e $V$ doit \^etre de la forme
$$n=p+m+pm$$
o\`u $p$ est le rang du $p$-syst\`eme de contact. En effet on a $q(p+1)=p(n-p)$ soit $n-p=q+q/p$. En posant $q=mp$, on obtient bien $n-p=mp+m.$ Par exemple, si $n=3$, alors $p=q=1$, si $n=4$ il n'y a pas de $p$-syst\`eme de contact. Si $n=5$, alors $p=1$ et $q=2$, ou bien $p=2$ et $q=2$. 
\begin{theo}\cite{Go-Ha}
Soit $\pp$ un $p$-syst\`eme de contact sur la vari\' et\' e $V$ de dimension $n=p+m+mp$. Pour tout $x \in V$, il existe un ouvert $U$, appel\' e ouvert de Darboux, de coordonn\' ees locales $\{x_1,\cdots,x_p,y_1,\cdots,y_m,z_1,\cdots,z_{pm}\}$ tel que $\pp$ soit repr\' esent\' e dans $U$ par le syst\`eme de Pfaff
$$
\left\{
\begin{array}{l}
\alpha_1=dx_1 +z_1dy_1+\cdots+z_{(m-1)p+1}dy_p\\
\alpha_2=dx_2+z_2dy_1+\cdots+z_{(m-1)p+2}dy_p\\
\cdots \\
\alpha_p=dx_p+z_pdy_1+\cdots+z_{mp}dy_p.
\end{array}
\right.
$$
\end{theo}  
Consid\' erons le syst\`eme diff\' erentiel de rang $p$ donn\' e par les champs $X_i$ dont les expressions dans l'ouvert de Darboux $U$ sont $X_i=\frac{\delta}{\delta x_i}$ pour $i=1,\cdots,p$. Ces champs v\' erifient $i(X_i)d\alpha_j=0$ et $\alpha_i(X_i)=1$.  Notons par $\mathcal{R}$ la distribution int\' egrale d\' efinie par ces champs et par $P=\pp / \mathcal{R}$ le fibr\' e quotient. Si $\mathcal{C}^\infty_{\mathcal{R}}$ est l'alg\`ebre des fonctions diff\' erentiables $f$ sur $V$ v\' erifiant $X_i(f)=0$ pour $i=1,\cdots,p$ alors $(P,\mathcal{C}^\infty_{\mathcal{R}})$ est muni d'une structure d'alg\' ebro\"{\i}de de Lie.

\section{Version alg\' ebrique des alg\' ebro\"{\i}des de Lie: les alg\`ebres de Lie-Rinehart}

Lorsque nous avons pr\' esent\' e la structure d'alg\' ebro\"{\i}de de Lie du fibr\' e tangent \`a une vari\' et\' e, nous avons mis en \' evidence deux structures alg\' ebriques, une alg\`ebre de Lie et une alg\`ebre associative avec des relations entre elles. Nous allons formaliser alg\' ebriquement ces structures en se passant du mod\`ele g\' eom\' etrique. On parle ainsi des alg\`ebres de Palais \' egalement appel\' ees alg\`ebre de Lie-Rinehart.

\subsection{Alg\`ebres de Palais ou de Lie-Rinehart}

\begin{defi}
On appelle alg\`ebre de Palais ou alg\`ebre de Lie-Rinehart un couple de deux $\K$-alg\`ebres $(L,A)$ o\`u
\begin{enumerate}
\item $A$ est une $\K$-alg\`ebre associative
\item $L$ est une $\K$-alg\`ebre de Lie et il existe un morphisme d'alg\`ebre de Lie $\rho: L \rightarrow Der(A)$ dans l'alg\`ebre des d\'erivations de $A$
ce qui donne une action \`a gauche de $L$ sur $A$:
$$ (X,a) \in L \times A \rightarrow X\cdot a = \rho(X)(a)$$
\item $L$ est un $A$-module c'est-\`a-dire est donn\' ee une multiplication externe
$$(a,X) \in A \times L \rightarrow aX \in L$$
\item Ces deux actions v\' erifient les relations suivantes:
\begin{equation}\label{LR}
[X,aY]=\rho(X)(a)Y+a[X,Y],  \ \ \ 
\rho(aX)(b)=a(\rho(X)(b))
\end{equation}
pour tout $X,Y \in L$ et $a,b \in A$
\end{enumerate}
o\`u $ab$ d\'esigne le produit dans $A$ et $[X,Y]$ celui dans $L$.
\end{defi}
Notons qu'en g\' en\' eral, une alg\`ebre de Lie-Rinehart est une alg\`ebre nonassociative. Si on note par $\diamond$ la multiplication sur cette alg\`ebre, c'est-\`a-dire $X \diamond Y = [X,Y], a\diamond b = ab$, $X \diamond a=X \cdot a$ et $a \diamond X=aX$, alors l'associateur de $\diamond$:
$$\mathcal{A}_\diamond(x,y,z)=x \diamond(y \diamond z)-(x \diamond y) \diamond z$$
v\'erifie:
$$
\begin{array}{lll}
\mathcal{A}_\diamond(X,a,Y)&=&X \diamond(aY)-(Xa)\diamond Y=[X,aY]-(X\cdot a)Y=a[X,Y],\\
\mathcal{A}_\diamond(a,X,Y)&=&a[X,Y]-[aX,Y]=Y(a)X,\\
\mathcal{A}_\diamond(X,Y,a)&=&X(Y(a))-[X,Y](a)=-Y(X(a)).,\\
\mathcal{A}_\diamond(X,a,b)&=& aX(b),\\
\mathcal{A}_\diamond(a,X,b)&=&0,\\
\mathcal{A}_\diamond(a,b,X)&=&0
\end{array}
$$
et n'est en g\'en\'eral pas nul. La condition $\mathcal{A}_\diamond(X,a,b)= aX(b)$ montre qu'il ne peut exister non plus des relations de sym\'etries sur l'associateur.
\medskip

\noindent{\bf Exemples d'alg\`ebres de Lie-Rinehart.}
\begin{enumerate}
\item
Le couple $(\Gamma(V), \mathcal{C}^\infty (V))$ o\`u $V$ est une vari\' et\' e diff\' erentielle est une alg\`ebre de Lie-Rinehart.
\item Consid\' erons l'alg\`ebre associative de dimension $3$ donn\' ee dans une base $\{e_1,e_2,e_3\}$ par
\begin{equation}\label{ex1}e_1e_i=e_ie_1=e_i, i=1,2,3, \ \ e_3e_3=e_2
\end{equation}
les produits non d\' efinis \' etant nuls. Consid\' erons l'alg\`ebre de Lie $L=Der(A)$ et donc $\rho=Id$. Elle est de dimension $2$, engendr\' ee par $X,Y$ les d\' erivations donn\' ees par
$$X(e_1)=0,X(e_2)=2e_2,X(e_3)=e_3, \ Y(e_i)=0, i=1,2, \ Y(e_3)=e_2$$
ce qui donne, d'apr\`es (\ref{LR}):
$$e_1X=X,e_1Y=Y, \ e_2X=e_2Y=0,\ e_3X=Y,e_3Y=0.$$
On en d\' eduit l'alg\`ebre de Lie-Rinehart de dimension $5$ dont la table de multiplication est
$$
\begin{array}{|c|c|c|c|c|c|}
\hline
 & X & Y & e_1 & e_2 &e_3 \\
\hline
X &0&Y &0&2e_2 &e_3 \\
\hline
Y & -Y & 0 & 0&0&e_2 \\
\hline
e_1 & X & Y& e_1 & e_2 & e_3 \\
\hline
e_2 &0&0&e_2 &0&0\\
\hline
e_3 & Y &0& e_3 &0&e_2 \\
\hline
\end{array}
$$ 
\item Plus g\'en\'eralement, soit $A$ une alg\`ebre associatie commutative de dimension finie $n$, $L$ l'alg\`ebre des d\'erivations de $A$. Comme $A$ est commutative, pour tout $a,b \in A$ et $X \in L$, la condition $(aX)(b)=aX(b)$ implique que $aX$ est bien une d\'erivation de $A$. En effet
$$(aX)(bc)=a(X(bc))=a(X(b)c+bX(c))=(aX)(b)c+b(aX)(c).$$
Comme $L$ co\"{\i}ncide dans ce cas avec $Der(A)$, la d\'erivation $aX$ est bien dans $L$. Dans le cas g\'en\'eral, c'est une contrainte sur $L$.  Dans l'exemple pr\'ec\'edent, si nous prenons pour $L$ la sous-alg\`ebre engendr\'ee par $X$, la condition $e_3X \in L$ n'est pas r\'ealis\'ee. Ainsi, dans cet exemple, nous voyons que les seules, \`a isomorphisme pr\`es, alg\`ebres de Lie-Rinehart $(A,L)$ lors que $A$ est donn\'ee par (\ref{ex1}) correspondent \`a $L=Der(A)$ ou $L=\K\{Y\}$.

\end{enumerate}

\subsection{Cas o\`u l'action de  $L$ sur $A$ est fid\`ele}
Par d\'efinition d'une alg\`ebre de Lie-Rinehart, il existe une action de $L$ sur $A$ associ\'ee au morphisme $\rho$ d'alg\`ebres de Lie. Nous dirons que cette action est fid\`ele si $X \cdot a =0$ pour tout $a \in A$ implique $X=0$.  Comme $X.a=\rho(X)(a)$ et $\rho(X) \in Der(A)$, cette action est fid\`ele si $\rho$ est injective. Si l'action $\rho$ est fid\`ele, nous dirons que l'alg\`ebre de Lie-Rinehart correspondante est fid\`ele. C'est en particulier le cas dans les alg\' ebro\"{\i}des d\' efinis pr\' ec\' edemment. Dans la d\' efinition de la structure de Lie-Rinehart, la condition (4) pr\' ecise deux identit\' es:
$$[X,aY]=X(a)Y+a[X,Y]$$
$$(aX)(b)=a(X(b)).$$
\begin{lemm}
Supposons que l'alg\`ebre de Lie-Rinehart $(L,A)$ soit fid\`ele,  alors
 $$(aX)(b)-a(X(b))=0 \Rightarrow [X,aY]-X(a)Y-a[X,Y]=0.$$
\end{lemm}
\dd En effet
$$[X,aY]b=X(aY(b))-(aY)(X(b))=X(a)Y(b)+a(X(Y(b)))-(aY)(X(b)).$$
Comme par hypoth\`ese $(aX)(b)-a(X(b))=0$ pour tout $a,b \in A$ et $X \in L$, on obtient
$$[X,aY]b=X(a)Y(b)+a(X(Y(b)))-a(Y(X(b)))=X(a)Y(b)+a(X(Y(b))-Y(X(b))$$
c'est-\`a-dire
$$[X,aY]b=(X(a)Y)b+a[X,Y]b.$$
Ainsi
$(aX)(b)=a(X(b))$ implique  $([X,aY]-X(a)Y+a[X,Y])(b)=0$ pour tout $b \in A$. Mais par hypoth\`ese, l'action de  $L$ sur $A$ est suppos\' ee fid\`ele ce qui implique
$$[X,aY]-X(a)Y+a[X,Y]=0.$$

Nous pouvons ainsi voir une alg\`ebre de Lie-Rinehart fid\`ele comme une alg\`ebre, not\'ee $\mathcal{LA}$, nonassociative $\Gamma$-gradu\' ee o\`u $\Gamma$ est un magma associatif \`a deux \' el\' ements $\alpha,\beta$ chacun d'eux \' etant un \' el\' ement neutre \`a gauche, dont le produit v\' erifie
\begin{enumerate}
\item $uv=(-1)^{d(\gamma)}vu$ pour tout $u,v$ homog\`enes, i.e. $u,v \in \mathcal{LA}_\gamma$ with $d(\alpha)=-d(\beta)=1$,
\item $\mathcal{A}(u,v,w)=\rho_1(\gamma_1,\gamma_2,\gamma_3)(wu)v+\rho_2(\gamma_1,\gamma_2,\gamma_3)v(uw)$

avec
\begin{enumerate}
\item $\rho_1(\gamma_1,\gamma_2,\gamma_3)=\rho_2(\gamma_1,\gamma_2,\gamma_3)=0$ when $\gamma_i,\gamma_j \in \mathcal{LA}_{\alpha}$
\item $\rho_1(\gamma_1,\gamma_2,\gamma_3)=-1, \ \ \rho_2(\gamma_1,\gamma_2,\gamma_3)=0$ when $(\gamma_1,\gamma_2,\gamma_3)=(\alpha,\beta,\beta)$
\item $\rho_1(\gamma_1,\gamma_2,\gamma_3)=0, \ \ \rho_2(\gamma_1,\gamma_2,\gamma_3)=-1$ when $(\gamma_1,\gamma_2,\gamma_3)=(\beta,\alpha,\beta)$ or $(\beta,\beta,\alpha)$ 
or $(\beta,\beta,\beta)$
\end{enumerate}
\end{enumerate}

\subsection{Alg\`ebres de Poisson et de Lie-Rinehart}
Nous avons rappel\' e ci-dessus la notion de vari\' et\' e de Poisson et la structure d'alg\' ebro\"{\i}de de Lie qui lui est attach\' ee.  Nous allons rappeler maintenant la version alg\' ebrique, c'est-\`a-dire celle d'alg\`ebres de Poisson.
\begin{defi}
Une $\K$-alg\`ebre de Poisson est la donn\' ee d'un triplet $(A, \cdot, \{.,.\})$
o\`u
\begin{enumerate}
\item $(A,\cdot)$ est une alg\`ebre associative commutative (on notera pour simplifier le produit $x \cdot y$ par $xy$,
\item $(A,\{.,.\})$ est une alg\`ebre de Lie,
\item les deux multiplications sont reli\' ees par la formule de Leibniz:
$$\{x,yz\}=y\{x,z\}+\{x,y\}z$$
pour tout $x,y,z \in A.$
\end{enumerate}
\end{defi}
Soit $(A,\cdot,\{.,.\})$ une alg\`ebre de Poisson. Supposons que le centre de l'alg\`ebre de Lie $(A,\{,\})$ soit r\'eduit \`a z\'ero. Pour tout $x \in A$, soit $X_x$ l'endomorphisme lin\' eaire d\' efini par
$$X_x(y)=\{x,y\}$$
pour tout $y \in A$. L'identit\' e de Leibniz implique
$$X_x(yz)=X_x(y)z+yX_x(z)$$
et donc $X_x$ est une d\' erivation de l'alg\`ebre associative $(A,\cdot)$. Notons par $L$ l'ensemble des $X_x$ lorsque $x$ parcourt $A$. C'est un sous-espace vectoriel de $Der(A,\cdot)$. De plus, si $[.,.]$ d\' esigne le crochet de Lie de l'alg\`ebre $Der(A,\cdot)$, alors
$$[X_x,X_y]=X_{\{x,y\}}.$$
En effet, si $z \in A$, alors
$$[X_x,X_y](z)=X_x(X_y(z))-X_y(X_x(z))=\{x,\{y,z\}\}-\{y,\{x,z\}\}.$$
Comme $\{.,.\}$ v\' erifie l'identit\' e de Jacobi, on obtient
$$[X_x,X_y](z)=\{\{x,y\},z\}.$$
Ainsi $L$ est une sous-alg\`ebre de Lie de $Der(A,\cdot)$. Notons que si $A$ est de dimension finie, et si $\{e_1,\cdots,e_n\}$ en est une base, alors si $x=\sum x_ie_i \in A$, nous avons $X_x=\sum x_i X_i$ o\`u $X_i=X_{e_i}$.  Comme $X_x(y)=\{x,y\}$, l'application li\' eaire $X_x$ est nulle si et seulement si $x$ est dans le centre  $Z_P$ de l'alg\`ebre de Lie $(A,\{.,.\})$ et $\dim L=\dim A -\dim Z_P.$

Consid\' erons \`a pr\' esent la paire $(L,A)$. On a bien:
\begin{enumerate}
\item $(A,\cdot)$ est une alg\`ebre associative commutative,
\item $L$ est une alg\`ebre de Lie, sous alg\`ebre de l'alg\`ebre des d\' erivations de $(A,\cdot)$.
\item Il existe une action de $L$ sur $A$  donn\' ee par
$$X_ x y= X_x(y)=\{x,y\}.$$
\end{enumerate}
D\' efinissons une action de $A$ sur $L$. Soient $y \in A$ et $X_x \in L$, alors $yX_x$ est l'\' el\' ement de $L$ d\' efini par
$$(yX_x)(z)=y(X_x(z))=y\{x,z\}.$$On en d\' eduit
$$
\begin{array}{ll}
[X_x,aX_y]z &= X_x(aX_y(z))-aX_y(X_x(z)\\
&= X_x(a)X_y(z)+aX_x(X_y(z))-aX_y(X_x(z))\\
&=X_x(a)X_y(z)+a[X_x,X_y](z)
\end{array}
$$ 
et donc
$$([X_x,aX_y]-X_x(a)X_y(z)+a[X_x,X_y])(z)=0$$
pour tout $z \in A.$ Mais pour tout $u,v \in A$, $X_uv=X_u(v)=\{u,v\}$ et $X_u(v)=0$ pour tout $v$ implique que $u$ soit dans le centre de l'alg\`ebre de Lie $(A,\{,\})$. D'apr\`es l'hypoth\`ese, on en d\'eduit que  l'application lin\' eaire $X_u$ est nulle. Ainsi
$$[X_x,aX_y] - X_x(a)X_y-a[X_x,X_y]=0$$
et nous r\' ecup\' erons une alg\`ebre de Lie-Rinehart. Bien entendu la r\' eciproque n'est pas vraie sans conditions en particulier sur les dimensions des alg\`ebres formant le couple d'alg\`ebres de Lie-Rinehart.

\subsection{Classification en petite dimension}
Pour \'elaborer une classification en dimension finie des alg\`ebres de Lie-Rinehart sur $\K=\C$ ou sur un corps de caract\'eristique $0$ et alg\'ebriquement clos, nous pouvons nous baser sur les listes des alg\`ebres associatives commutatives connues et d\'eterminer, non pas l'alg\`ebre de Lie $L$ mais son image par $\rho$, comme le noyau de $\rho$ est un id\'eal de $L$, la d\'etermination de $L$ s'en suivra. Dans ce qui suit, $L$ d\'esignera donc en fait son image par $\rho$.
\subsubsection{$\dim A =2$}
Soit $A$ une alg\`ebre associative commutative de dimension $2$ et notons par $\{e_1,e_2\}$ une base de cette alg\`ebre. 
\begin{enumerate}
\item $$e_1e_i=e_ie_1=e_i, \ \ e_2e_2=e_2.$$
 Dans ce cas les d\' erivations $X_i=X_{e_i}$ sont toutes nulles et l'alg\`ebre de Lie-Rinehart se r\' eduit \`a la seule alg\`ebre associative $A$. 
\item  $$e_1e_i=e_ie_1=e_i, \ \ e_2e_2=0.$$
Dans ce cas $Der(A)$ est de dimension $1$ et engendr\' ee par la d\' erivation $X$ donn\' ee par
$$X(e_1)=0, \ X(e_2)=e_2.$$
Prenons pour $L$ l'alg\`ebre de Lie unidimensionelle engendr\' ee par $X$. Il nous reste \`a d\' efinir $e_1X$ et $e_2X$. Le calcul des images des vecteurs de base montre que
$$e_1X=X, \ \ e_2X=0.$$
On a donc la structure de l'alg\`ebre de Lie-Rinehart d\' ecrite dans le tableau suivant:
$$\begin{array}{c||c||c|c}
 &X & e_1 & e_2 \\
\hline
X & 0 & 0 & e_2 \\
\hline
e_1 & X & e_1 & e_2 \\
\hline
e_2 &0 & e_2 & 0
\end{array}
$$
\item $$e_1e_1=e_2.$$
Rappelons que les produits non d\' efinis sont suppos\' es nuls.  L' alg\`ebre des d\'erivations de $A$ est de dimension $2$ et engendr\' ee par les deux d\' erivations
$$X(e_1)=e_1,X(e_2)=2e_2, \ \ Y(e_1)=e_2, Y(e_2)=0.$$
Comme $e_1X=Y$, on en d\' eduit les structures de Lie-Rinehart suivantes:
$$\begin{array}{c||c|c||c|c}
 &X  & Y &e_1 & e_2 \\
\hline
X & 0 & Y & e_1 & 2e_2 \\
\hline
Y & 0 & 0 & e_2 & 0 \\
\hline
e_1 & Y & 0& e_2 &0 \\
\hline
e_2 &0& 0 &0 & 0
\end{array}, \ \ \ \ \ 
\begin{array}{c||c||c|c}
 & Y &e_1 & e_2 \\
\hline
Y& 0 & e_2 & 0 \\
\hline
e_1 &  0& e_2 &0 \\
\hline
e_2 & 0 &0 & 0
\end{array},
$$
\item $$e_1e_1=e_1.$$
Comme $Der(A)$ est une alg\`ebre de Lie est de dimension $1$  engendr\' ee par la d\' erivation
$$X(e_1)=0,X(e_2)=e_2,$$
on en d\' eduit la structure de Lie-Rinehart suivante:
$$\begin{array}{c||c||c|c|}
 &X   &e_1 & e_2 \\
\hline
X & 0  & 0 & e_2 \\
\hline
e_1 &  0& e_1 &0 \\
\hline
e_2 &0 &0 & 0
\end{array}
$$
\item $$e_ie_j=0.$$
L'alg\`ebre des dérivations est l'algèbre de Lie $gl(2)$ des matrices d'ordre $2$. Consid\'erons la base canonique $\{X_1,X_2,X_3,X_4\}$ de cette algèbre de Lie. Comme $(e_iX)e_j=e_i(X(e_j)=0$, tous les vecteurs $e_iX_j$ sont nuls. L'algèbre de Lie-Rinehart pleine (voir le paragraphe qui suit) associée est donc
$$\begin{array}{c||c|c|c|c||c|c|}
 &X_1  & X_2 & X_3 & X_4 &e_1 & e_2 \\
\hline
X_1 & 0 & X_2 & -X_3 &0 & e_1 & 0 \\
\hline
X_2 & -X_2 & 0 & X_1-X_4 &X_2 & 0 & e_1 \\
\hline
X_3 & X_3 & -X_1+X_4 & 0 &-X_3 & e_2 & 0 \\
\hline
X_4 & 0 & -X_2 & X_3 &0 & 0 & e_2 \\
\hline
e_1 &0 & 0&0&0& 0 &0 \\
\hline
e_2 &0& 0&0&0 &0 & 0
\end{array}
$$
Pour toute sous-alg\`ebre $L'$ de $L$, le couple $(L',A)$ est aussi une algèbre de Lie-Rinehart. En particulier l'algèbre engendrée par $\{X_3,X_4,e_1,e_2\}$ est une algèbre de Lie-Rinehart attachée à l'algèbre de Poisson définie par
$$e_ie_j=0, \ \{e_1,e_2\}=e_2$$
car nous avons dans ce cas
$$X_{e_1}(e_1)=0,X_{e_1}(e_2)=e_2, \ \ X_{e_2}(e_1)-=e_2, X_{e_2}(e_2)=0$$
et $[X_{e_1},X_{e_2}]=X_{e_2}$.
On en d\' eduit la structure de Lie-Rinehart associ\' ee \`a l'alg\`ebre de Poisson:
$$\begin{array}{c||c|c||c|c|}
 &X  & Y &e_1 & e_2 \\
\hline
X & 0 & Y & 0 & e_2 \\
\hline
Y & -Y & 0 & -e_2 & 0 \\
\hline
e_1 &0 & 0& 0 &0 \\
\hline
e_2 &0& 0 &0 & 0
\end{array}
$$

\end{enumerate}

\subsubsection{Dimension 3} La liste des alg\`ebres de Poisson de dimension $3$ est donn\' ee dans \cite{Re1, G.R.Poiss}. Nous n'allons pas passer en revue toute cette liste, elle est trop longue, mais regarder la première de la liste. On considère l'algèbre associative commutative de dimension $3$ donnée par
$$e_1e_i=e_i, \ i=1,2,3. $$

L'algèbre de dérivation est de dimension $4$. Elle est engendrée par les dérivations $X_i$, $i=1,\cdots 4$, donnée par
$$X_1(e_2)=e_2,X_2(e_2)=e_3,X_3(e_2)=e_2,X_4(e_3)=e_3$$
les images non définies étant nulles. On en déduit que $e_iX_j=0$ dès que $i =2,3$ et $e_1X_i=X_i$, d'où l'algèbre
$$\begin{array}{c||c|c|c|c||c|c|c}
 &X_1&   X_2 & X_3 & X_4 &e_1 & e_2 & e_3 \\
\hline
X_1 & 0 & -X_2  & X_3&0 & 0 & e_2 &0 \\
\hline
X_2 & X_2 &0 & X_4-X_1  & -X_2&0 & e_3 &0 \\
\hline
X_3 & -X_3 &X_1-X_4 & 0  & X_3&0 & 0 &e_2 \\
\hline
X_4 & 0 &X_2 & -X_3  & 0&0 & 0 &e_3 \\
\hline
e_1 &X_1&   X_2 & X_3 & X_4& e_1 &e_2& e_3\\
\hline
e_2 &0&   0 & 0 & 0& e_2 &0& 0\\
\hline
e_3 &0&   0 & 0 & 0& e_3 &0& 0\\
\end{array}
$$
\noindent Remarque. Consid\'erons l'algèbre de Poisson construite sur $A$:
$$e_1e_i=e_i, \ i=1,2,3, \ \ \ \{e_2,e_3\}=e_3.$$
Dans ce cas $X_{e_1}=0$ et $X_{e_2}(e_3)=e_3$, $X_{e_3}(e_2)=-e_3$ ces d\' erivations \' etant nulles sur les autres vecteurs de base soit
$$X_{e_2}=X_4, X_{e_3}=-X_2.$$
On en d\' eduit l'alg\`ebre de Lie-Rinehart attach\' ee \`a cette structure de Poisson
$$\begin{array}{c||c|c||c|c|c|}
 &X  & Y &e_1 & e_2 & e_3 \\
\hline
X & 0 & Y & 0 & 0 &e_3 \\
\hline
Y & -Y & 0 &0&-e_3 &0 \\
\hline
e_1 &X & Y& e_1 &e_2& e_3\\
\hline
e_2 &0& 0 &e_2 & 0 & 0\\
\hline
e_3&0& 0 &e_3& 0 &0 \\
\end{array}
$$

\section{Alg\`ebres de Lie-Rinehart pleines}

L'examen des exemples pr\' ec\' edents montre l'int\' er\^et d'\' etudier, au moins en dimension finie, les alg\`ebres de Lie-Rinehart $(A,L)$ lorsque $L$ co\"{\i}ncide avec l'alg\`ebre de Lie des d\' erivations de l'alg\`ebre associative $A$. Ceci conduit \`a la d\' efinition suivante
\begin{defi}
On dit qu'une alg\`ebre de Lie-Rinehart est pleine lorsque l'alg\`ebre de Lie $L$ est l'alg\`ebre des d\' erivations de $A$.
\end{defi}
Dans ce cas nous avons
\begin{propo}
Soit $A$ une $\K$-alg\`ebre associative commutative de dimension finie. Il existe, \`a isomorphisme pr\`es, une et une seule alg\`ebre de Lie-Rinehart pleine dont la partie associative est $A$.
\end{propo}
En effet $Der(A)$ est enti\`erement d\' efinie par $A$ ainsi que la structure du $L$-module $A$.  Enfin les \' el\' ements de $L$ du type $aX$ avec $a \in A$ et $X \in L$ sont d\' etermin\' es par leurs images des vecteurs d'une base $\{e_i\}$ de $A$ par la formule
$$(aX)(e_i)=a(X(e_i)).$$
Il est clair que cette alg\`ebre de Lie-Rinehart est aussi fid\`ele. Dans les exemples pr\' ec\' edents nous avons d\' etermin\' e toutes ces alg\`ebres lorsque $\dim A=2$. 

Un exemple intéressant est celui où l'alg\`ebre associative $A$ est l'alg\`ebre triviale $O_n$ et donc $Der(A)=gl(n)$. L'alg\`ebre de Lie-Rinehart concerne le couple $(O_n,gl(n))$. Notons par $\mu$ la multiplication d\' efinie sur $A \oplus L$ \`a partir de la structure de Lie-Rinehart de $(A,L)$.
\begin{propo}
Soit $(O_n \oplus gl(n),\mu)=(O_n,gl(n))$ l'alg\`ebre de Lie-Rinehart pleine avec $A=O_n$. Alors si $\psi_\mu$ est l'application anti-sym\' etrique d\' efinie par
$$\psi_\mu(a,X)=-\psi_\mu(X,a)=-X(a), \ \psi_\mu(X,Y)=\mu(X,Y), \ \ \psi_\mu(a,b)=0$$
pour tout $a,b \in O_n$ et $X,Y \in gl(n)$, l'alg\`ebre $(O_n\oplus gl(n),\psi_\mu)$ est une alg\`ebre de Lie.
\end{propo}   
\dd En effet, la condition de Jacobi sur $\psi_\mu$ est satisfaite si et seulement si elle est satisfaite sur les triplets du type $(X,Y,a)$.
$$\psi_\mu(\psi_\mu(X,Y),a)+\psi_\mu(\psi_\mu(Y,a),X)+\psi_\mu(\psi_\mu(a,X),Y)=[X,Y](a)-X(Y(a))+Y(X(a)=0$$
d'o\`u la proposition.

On remarquera que $\psi_{\mu}$ ne co\"{\i}ncide pas avec l'application anti-sym\' etrique attach\' ee \`a $\mu$. La structure d'alg\`ebre sur $O_n \oplus gl(n)$ attach\' ee \`a cette multiplication antisym\' etrique est  une structure d'alg\`ebre de type Lie au sens de \cite{ma}.

\medskip

Pour une alg\`ebre de Lie-Rinehart g\' en\' erale, cette application antisym\' etrique $\psi$ est donn\' ee par
$$
\left\{
\begin{array}{l}
 \psi(X,a)=-\psi(a,X)=X(a)-aX,\\
\psi(X,Y)=[X,Y],\\
\psi(a,b)=0
\end{array}
\right.
$$
pour tout $a,b \in A$ and $X,Y \in L$.  Si $J_\psi$ d\' esigne le Jacobian de $\psi$ (c'est-\`a-dire l'identit\' e de Jacobi relative \`a $\psi$), on a 
$$
\left\{
\begin{array}{l}
\medskip
 J_\psi(a,b,c)=0,\\
\medskip
 J_\psi(a,b,X)=-J_\psi(a,X,b)=J_\psi(X,a,b)=bX(a)-aX(b),\\
\medskip
J_\psi(a,X,Y)=-J_\psi(X,a,Y)= J_\psi(X,Y,a)= a[X,Y].
\end{array}
\right.
$$
Consid\' erons l'application trilin\' eaire anti-sym\' etrique $\Phi$ d\' efinie par
$$
\left\{
\begin{array}{l}
\Phi (a,b,c)=0,\\
\Phi (a,b,X)= bX(a)-aX(b),\\
\Phi (a,X,Y)=a[X,Y],\\
\Phi (X,Y,Z)= 0
\end{array}
\right.
$$
Bien entendu, on a $J_\psi=\Phi.$ D\' efinissons $\Phi$ comme la d\' eriv\' ee d'une application bilin\' eaire. Pour cela consid\' erons l'application bilin\' eaire
$$
\left\{
\begin{array}{l}
\phi (a,b)=0,\\
\phi (a,X)= aX-X(a),\\
\phi (X,Y)=X \wedge Y\\
\end{array}
\right.
$$
pour tout $a,b \in A$ et $X,Y \in L$. Ici $ X\wedge Y$ d\' esigne le bivecteur agissant sur $A$ de la fa\c{c}on suivante:
$$X \wedge Y (a)=X(a)Y-Y(a)X.$$
Soit $\delta$ l'application qui \`a $\phi$ fait correspondre l'application trilin\' eaire
$$
\left\{
\begin{array}{ll}
\delta(\phi) (a,b,c)= &0,\\
\delta(\phi (a,b,X)= &\phi(a,b)X+\phi(b,X)a+\phi(X,a)b+\phi(ab-ba,X)+\phi(bX-Xb,a)\\
&+\phi(Xa-aX,b)\\
\delta(\phi (a,X,Y)=&\phi(a,X)Y+\phi(X,Y)a+\phi(Y,a)X+\phi(aX-Xa,Y)+\phi(X\wedge Y,a)\\
&+\phi(Ya-aY,X)\\
\delta(\phi (X,Y,Z)=&\phi(X,Y)Z+\phi(Y,Z)X+\phi(Z,X)Y+\phi(X \wedge Y,Z)+\phi(Y\wedge Z,X)\\&+\phi(Z \wedge X,Y).\\
\end{array}
\right.
$$
On a donc 
$$
\left\{
\begin{array}{ll}
\delta(\phi) (a,b,c)= &0,\\
\delta(\phi (a,b,X)= & bX(a)-aX(b)=\Phi (a,b,X)\\
\delta(\phi (a,X,Y)=&a[X,Y]=\Phi(a,X,Y)\\
\delta(\phi (X,Y,Z)=& 0).\\
\end{array}
\right.
$$
Ainsi l'application $\psi$ satisfait l'identit\' e de Jacobi \`a un "cobord" pr\`es. Nous d\' evelopperons plus tard cette approche vers une structure de $L_\infty$-alg\`ebre.

\section{Alg\' ebres de Courant}

Nous n'allons pas développer ici cette extension de notion d'algèbre de Lie. Intéressons nous à l'exemple suivant: soit $\g$ une algèbre de Lie de dimension finie et soit $\g^*$ le dual vectoriel assimilé à l'espace des formes invariantes à gauche sur un groupe de Lie associé.  Soit $\frak{a} =\g \oplus \g^\star$ munie de l'application bilinéaire:
$$[[(X_1,\omega_1),(X_2,\omega_2)]]=([X_1,X_2],L(X_1)\omega_2)$$
o\`u $L(X)$ d\' esigne la d\' eriv\' ee de Lie des formes diff\' erentielles associ\' ees \`a $X$, soit $L(X)\omega=i(X)d\omega+d(i(X)\omega)$, $i(X)$ \' etant le produit int\' erieur. Comme $d(i(X)\omega)=0$, alors $L(X)\omega=i(X)d\omega$. Prenons par exemple pour $\g$ l'alg\`ebre de Lie r\' esoluble non ab\' elienne de dimension $2$. Ses constantes de structure sont donn\' ees par $[e_1,e_2]=e_2$ et donc $d\omega_1=0, d\omega_2=-\omega_1 \wedge \omega_2$. Posons  $v_1=(e_1,0), v_2=(e_2,0),v_3=(0,\omega_1),v_4=(0,\omega_2).$  On a donc la table de multiplication suivante:
$$\begin{array}{c||c|c|c|c}
 &v_1  & v_2 &v_3&v_4 \\
\hline
v_1 & 0 & v_2& 0 &-v_4 \\
\hline
v_2 & -v_2 & 0 & 0 & v_3 \\
\hline
v_3 & 0 & 0& 0 &0 \\
\hline
v_4 &0& 0 &0 & 0
\end{array}
$$
Ce crochet n'est pas un crochet de Lie.  Il vérifie l'identité
$$[u_i,[u_j,u_k]]-[[u_i,u_j],u_k]-[u_j,[u_i,u_k]=0$$
el l'algèbre correspondante est appelée algèbre de Leibniz ou algèbre de Loday. Une telle structure algèbrique est appelée une algèbre de Courant, modèle algébrique d'un algébro\"{\i}de de Courant. Rappelons-en la définition: Une algèbre de Courant sur une algèbre de Lie $\g$ est une algèbre de Leibniz $(L, [[,]])$
munie d'une application linéaire $\pi: L \rightarrow \g$ telle que
$$\pi[[u_1,u_2]]=[\pi(u_1),\pi(u_2]$$
pour tout $u_1,u_2 \in L.$
Elle est dite complète lorsque $\pi$ est surjective. 

Considérons l'algèbre de Courant $\g \oplus \g^*$, nous pouvons considérer $\g^*$ comme une algèbre associative abélienne. L'algèbre de Courant $\g \oplus \g^*$ s'identifie alors au couple $(\g,\g^*)$ formé d'une algèbre de Lie et d'une algèbre associative abélienne pour lequel on a une action triviale de $\g *$ sur $\g$ car 
$$[[(0, \omega), (X,0)]]=0,$$
d'une action de $\g$ sur $\g^*$:
$$[[(0,\omega)]]=L(X)\omega=i(X)d\omega$$
qui fait de $\g^*$ un $\g$-module. Ainsi $(\g,\g^*)$ est une algèbre de Lie-Rinehart. Dans ce cas, la multiplication $\diamond$ associée définie au paragraphe 3.1 est une multiplication d'algèbre de Leibniz.

\end{document}